\documentclass{article}

\usepackage{fullpage}
\usepackage{url}
\usepackage[utf8]{inputenc}
\usepackage{amsthm}
\usepackage{amsmath}
\usepackage{amssymb}
\usepackage{bbm}
\usepackage{physics}
\usepackage{mathtools}

\DeclarePairedDelimiter{\floor}{\lfloor}{\rfloor}

\newtheorem{thm}{Theorem}[section]
\newtheorem{cor}[thm]{Corollary}

\title{Additive energy, discrepancy and Poissonian $k$-level correlation}
\author{
  Guy Lachman\\
  \texttt{guylachman@mail.tau.ac.il}
  \and
  Shvo Regavim\\
  \texttt{shvoregavim@mail.tau.ac.il}
}
\date{\today}

\begin{document}

\maketitle

\begin{abstract}
    $k$-level correlation is a local statistic of sequences modulo 1, describing the local spacings of $k$-tuples of elements. For $k = 2$ this is also known as pair correlation. We show that there exists a well spaced increasing sequence of reals with additive energy of order $N^3$ and Poissonian $k$-level correlation for all integers $k \geq 2$, answering in the affirmative a question raised by Aistleitner, El-Baz, and Munsch. The construction is probabilistic, and so we do not obtain a specific sequence satisfying this condition. To prove this, we show that random perturbations of a sequence with small discrepancy gives, almost surely, a sequence with Poissonian $k$-level correlation, a fact which may be of independent interest.
\end{abstract}

\section{Introduction}
Throughout, we will let $\{ x \}, \norm{x}$ denote the fractional part of $x$ and the distance from $x$ to the nearest integer respectively. A sequence of reals $\left( x_n \right)_{n \geq 1}$ is \emph{equidistributed modulo 1} if for all $0 < a < b < 1$ we have
$$\lim_{N \to \infty} \frac{1}{N} \sum_{n = 1}^{N} \mathbbm{1} \left( \{ x_n \} \in [a, b] \right) = b - a$$
that is the number of $x_n$'s which lie in a certain interval is asymptotically proportional to the length of that interval. A sequence of reals $\left( x_n \right)_{n \geq 1}$ is said to have \emph{Poissonian pair correlation} if for all $s > 0$ we have
$$\lim_{N \to \infty} \frac{1}{N} \sum_{1 \leq m \neq n \leq N} \mathbbm{1} \left( \norm{x_m - x_n} < \frac{s}{N} \right) = 2 s$$
which essentially means that the local spacing at scale $\frac{1}{N}$ between pairs of elements is the same as in the random case. Poissonian pair correlation implies equidistribution (see \cite{pair->equi1}, \cite{pair->equi2}), and so it can be thought as a finer, more local statistic.

Given a sequence of reals $x_n$, a natural number $N > 0$ and a real $\gamma > 0$, define the additive energy at scale $\gamma$ to be
$$E_{N, \gamma} = \# \left\{ 1 \leq a, b, c, d \leq N \ \big\lvert \ \abs{x_a + x_b - x_c - x_d} < \gamma \right\}$$
and when the scale is not written we let $E_{N} = E_{N, 1}$. We always have $E_{N, \gamma} \geq N^2$ by the contribution of the diagonal $a = c, b = d$. If for example we are given a spacing condition like $x_{n + 1} - x_n \geq \delta$ then $E_{N, \gamma} \leq \frac{2 \gamma N^3}{\delta} + N^3$ as $a, b, c$ define $d$ up to $\frac{2 \gamma}{\delta} + 1$ options.

We say that an increasing sequence of real numbers $\left( x_n \right)_{n \geq 1}$ is \emph{well spaced} if $x_{n + 1} - x_n \geq 1$. Given a property $T$ of sequences, we say that a sequence $\left( x_n \right)_{n \geq 1}$ has the property \emph{metric $T$} if for almost all $\alpha \in \mathbb{R}$ in the sense of Lebesgue measure, the sequence $\left( \alpha x_n \right)_{n \geq 1}$ has the property $T$.

Aistleitner, El-Baz and Munsch in \cite{2009.08184} raised the problem of whether there exists an increasing sequence of reals $\left( x_n \right)_{n \geq 1}$ with additive energy of high order, that is $E_{N} \gg N^3$, and metric Poissonian pair correlation. Recently, Lutsko, Sourmelidis and Technau \cite{very small theta} have shown that for all $0 < \theta < \frac{1}{3}$ and for all $\alpha \neq 0$ the sequence $\alpha n^{\theta}$ has Poissonian pair correlation, and Rudnick and Technau \cite{small theta} have shown that for all $0 < \theta < 1$, the sequence $\left( n^{\theta} \right)_{n \geq 1}$ has metric Poissonian pair correlation, answering this question in the affirmative. However, as opposed to most of the sequences discussed in \cite{2009.08184}, these sequences are not well spaced, that is $x_{n + 1} - x_n \to 0$. We will show that there exist well spaced sequences, with additive energy $E_{N} \gg N^3$ and metric Poissonian pair correlation.

In \cite{2009.08184} it is conjectured that an increasing sequence of reals with additive energy $E_{N, \gamma} \gg \gamma N^4$ where $\frac{1}{N} \leq \gamma = \gamma (N) \leq 1$ is some function of $N$ is sufficient to ensure lack of Poissonian pair correlation. In the case of well spaced sequences, then for all $\gamma \leq 1$ we have $E_{N, \gamma} \ll N^3$, and so this criterion can only possibly hold when $\gamma \ll \frac{1}{N}$. The analogous question one can ask is given a well spaced sequence such that $E_{N, \gamma} \gg N^3$, how small must $\gamma$ be to ensure that $x_n$ does not have metric Poissonian pair correlation. More specifically, we ask
\\ \\
\textbf{Open Problem 1:} Show that a well spaced sequence with additive energy $E_{N, \gamma} \gg N^3$ does not have metric Poissonian pair correlation, where $\gamma = \frac{1}{N}$.
\\ \\
In the case where $x_n$ is an integer valued sequence, notice that $E_{N, \gamma} = E_N$ for all $\gamma \leq 1$. In an appendix to a paper by Aistleitner, Larcher, and Lewko \cite{bourgain}, Bourgain showed that if $x_n$ is an increasing sequence of integers with $E_N \gg N^3$ for infinitely many $N$, then $x_n$ does not have metric Poissonian pair correlation. This argument was later sharpened by Lachmann and Technau \cite{lachmann}, and essentially optimized by Larcher and Stockinger \cite{bourgain+}, who showed that in this case, there is no $\alpha$ such that the sequence $\alpha x_n$ has Poissonian pair correlation. We can adapt the argument of Bourgain to the real valued case, by defining $y_n = \frac{1}{C N} \floor{C N x_n}$ for some large constant $C$. Then, the large additive energy of $x_n$ with respect to a small scale transfers to many equalities of the form $y_a + y_b = y_c + y_d$, which then allows us to use the Balog-Szemeredi-Gowers theorem. This follows through to show that for all $c_1 > 0$ there is a constant $c_2 > 0$ such that if $E_{N, \frac{c_2}{N}} \geq c_1 N^3$ then the sequence $x_n$ does not have metric Poissonian pair correlation, but we have not managed to quite push it all the way to $E_{N, \frac{1}{N}} \gg N^3$.

In the other direction, we show the following:

\begin{thm}\label{real-result}

For each constant $c > 0$, there exists a well spaced sequence $x_n$ with additive energy $E_{N, \gamma} \gg N^3$ and metric Poissonian pair correlation, where $\gamma = \frac{\log N \left( \log \log N \right)^{1 + c}}{N}$.

\end{thm}

Our construction is probabilistic. Recall that the discrepancy of a sequence $x_n$ of numbers modulo 1 is defined as
$$D_N = \sup_{0 \leq a \leq b \leq 1} \abs{\frac{\# \left\{ n : x_n \in \left[ a, b \right] \right\}}{N} - \left( b - a \right)}$$
Discrepancy is a quantitative version of equidistribution: the more evenly distributed a sequence is modulo 1, the smaller the discrepancy.

Define $M_N = \max_{n \leq N} n D_n$. We will prove the following statement:

\begin{thm}\label{disc_pois}

Let $x_n$ be a sequence of numbers modulo 1, and let $g$ be a monotonically decreasing function, such that
$$\lim_{N \to \infty} \frac{g(N)}{D_N} = \infty$$
Furthermore, assume $g$ satisfies the following regularity conditions:
$$\lim_{N \to \infty} N g(N) = \infty$$
and
$$\lim_{N \to \infty} \frac{g \left( N \left( 1 + \frac{M_N}{N g(N)} \right) \right)}{g(N)} = 1$$
Let $z_1, z_2, \dots$ be independent random variables such that $z_i \sim \mathrm{Unif} \left[ -g(i), g(i) \right]$. Then, almost surely, the sequence $x_n + z_n$ has Poissonian pair correlation.

\end{thm}

Let us just mention that the regularity conditions on $g$ are just a technicality, and in all reasonable cases they are satisfied. A result of Schmidt \cite{lowdisc} states that there is some constant $c > 0$, such that infinitely often we have $D_N \geq \frac{c \log N}{N}$, and therefore it is reasonable to assume that $g(N) \geq \frac{\log N}{N}$ which is stronger than the first regularity condition. As for the second condition: intuitively, $\frac{M_N}{N}$ should be about $D_N$, and then $\frac{D_N}{g(N)} \to 0$, and so all we require is that when we multiply $N$ by a number which is $1 + o(1)$, the value of the function $g$ should change by $1 + o(1)$, which happens for most reasonable options for $g$. This regularity condition can be slightly weakened further, but for the reasons stated we have not attempted to optimize it.

This theorem suggests that for a "typical" sequence, Poissonian pair correlation (or other local statistics) are in some sense a "smoothed out" version of small discrepancy, and it is evident from the proof as well: the small discrepancy ensures that we are evenly distributed on a relatively small scale, and the random shifts by $z_n$ smooth out the sequence to give us pseudorandom behaviour at a scale of $\frac{1}{N}$, that is Poissonian pair correlation. The connection between these two statistics has been studied for example by Steinerberger \cite{stein} who has shown that a sequence that has a certain uniform version of Poissonian pair correlation has small discrepancy. Interestingly, Larcher and Stockinger \cite{neg} have established that many classical low discrepancy sequences do not have Poissonian pair correlation.

Let us show the deduction of Theorem \ref{real-result} from Theorem \ref{disc_pois}: Let $g(N) = \frac{\log N \left( \log \log N \right)^{1 + c}}{N}$. We will show that almost surely, a sequence of the form $x_n = 2 n + z_n$, where $z_n \sim \mathrm{Unif} \left[ -g(N), g(N) \right]$ are independent random variables, has metric Poissonian pair correlation. Clearly, this sequence satisfies $x_{n + 1} - x_n \geq 1$. This sequence has additive energy $E_{N, \gamma} \gg N^3$, where $\gamma = 10 g(N)$, because there are $\gg N^3$ quadruples $\frac{N}{2} \leq a, b, c, d \leq N$ such that $a + b = c + d$, and for each such quadruple we have $\abs{x_a + x_b - x_c - x_d} \leq 10 g(N)$. Now, by Fubini's theorem, in order to show that almost surely, for almost all $\alpha$ the sequence $\alpha x_n$ has Poissonian pair correlation, it is sufficient to show that for almost all $\alpha$, almost surely the sequence $\alpha x_n$ has Poissonian pair correlation. By a theorem of Beck \cite{beck} (though this specific case was proven first by Khintchine \cite{khintchine}), for almost all $\alpha$, the sequence $\alpha n$ has discrepancy
$$D_N = \mathcal{O} \left( \frac{\log N \left( \log \log N \right)^{1 + c/2}}{N} \right)$$
and a direct application of Theorem \ref{disc_pois} gives us the desired result.

In \cite{beck} Beck showed that for almost all $\alpha$, for infinitely many $N$ we have that the discrepancy of the sequence $\alpha n$ is at least $D_N \geq \frac{\log N \log \log N}{N}$, and so Theorem \ref{real-result} is the best we can get from Theorem \ref{disc_pois} with a sequence of the form $x_n = 2 n + z_n$. This raises a natural question, which to the best of our knowledge seems not to have been discussed in the literature:
\\ \\
\textbf{Open Problem 2:} Is there an increasing sequence of integers $x_n$ such that for almost all $\alpha$, the discrepancy of $\alpha x_n$ is $D_N \ll \frac{\log N \log \log N}{N}$?
\\ \\
It seems reasonable to conjecture that the answer to this question is negative: for example, Aistleitner and Larcher \cite{disc-ener} have shown that an increasing sequence of integers $x_n$ such that for almost all $\alpha, \ \alpha x_n$ has small discrepancy, must have large additive energy, which is a significant structural condition on the sequence. As we have stated above, Schmidt \cite{lowdisc} has shown that for any sequence modulo 1, for infinitely many $N$, the discrepancy is at least $D_N \geq \frac{\log N}{N}$, and so the best improvement to Theorem \ref{real-result} that our method involving Theorem \ref{disc_pois} could theoretically yield is $\gamma = f(N) \frac{\log N}{N}$, where $f$ is a slowly growing function which tends to infinity with $N$. Therefore, we ask
\\ \\
\textbf{Open Problem 3:} Is there a well spaced sequence $x_n$ with metric Poissonian pair correlation, with additive energy $E_{N, \gamma} \gg N^3$, where $\gamma \ll \frac{\log N}{N}$?
\\ \\
The proof of Theorem \ref{disc_pois} breaks down in this case if we try to apply the same method, which hints maybe that such sequences (if they exist) are atypical in some sense. We show a certain converse theorem to \ref{disc_pois}, which is not tight, but at least shows that our method in Theorem \ref{real-result} cannot get metric Poissonian pair correlation and additive energy at a scale of $\frac{\sqrt{\log N}}{N}$.

\begin{thm}\label{converse}

Let $0 < c \leq \frac{1}{2}$ be a real number and let $g(n) = \frac{\left( \log n \right)^c}{n}$. Let $z_1, z_2, \dots$ be independent random variables such that $z_i \sim \mathrm{Unif} \left[ -g(i), g(i) \right]$. Then, almost surely, the sequence $x_n = n + z_n$ does not have Poissonian pair correlation.

\end{thm}

For an integer $k \geq 2$, $k$-level correlation is a statistic which intuitively measures the local spacing at scale $\frac{1}{N}$ between $k$-tuples of elements. More formally, let $\mathcal{X}_k = \mathcal{X}_k (N)$ denote the set of $k$-tuples $\left( a_1, \dots, a_k \right)$ of distinct integers $1 \leq a_i \leq N$. Given $s = \left( s_1, \dots, s_{2 k - 2} \right) \in \mathbb{R}^{2 k - 2}$ where $s_{2 i - 1} < s_{2 i}$ for all $1 \leq i < k$, let $\mathbbm{1}_{s, N} \left( y_1, \dots, y_{k - 1} \right)$ denote the indicator function of the event $\{ y_1 \} \in \left[ \frac{s_1}{N}, \frac{s_2}{N} \right], \dots, \{ y_{k - 1} \} \in \left[ \frac{s_{2 k - 3}}{N}, \frac{s_{2 k - 2}}{N} \right]$ where the intervals are taken modulo 1, that is for our purposes $\frac{3}{4} \in \left[ - \frac{1}{3}, \frac{1}{3} \right]$. A sequence of reals $\left( x_n \right)_{n \geq 1}$ is said to have \emph{Poissonian $k$-level correlation} if for all $s_1, s_2, \dots, s_{k - 1} > 0$ we have
$$\lim_{N \to \infty} \frac{1}{N} \sum_{a \in \mathcal{X}_k} \mathbbm{1}_{s, N} \left( x_{a_1} - x_{a_2}, x_{a_1} - x_{a_3}, \dots, x_{a_1} - x_{a_k} \right) = \left( s_2 - s_1 \right) \left( s_4 - s_3 \right) \cdots \left( s_{2 k - 2} - s_{2 k - 3} \right)$$
Analogously to pair correlation, Poissonian $k$-level correlation means that the local spacing at scale $\frac{1}{N}$ between $k$-tuples of elements is the same as in the random case. Usually in the literature Poissonian $k$-level correlation is defined by the property that for any smooth, compactly supported $f \in C_{c}^{\infty} \left( \mathbb{R}^{k - 1} \right)$ we have
$$\lim_{N \to \infty} \frac{1}{N} \sum_{a \in \mathcal{X}_k} \sum_{m \in \mathbb{Z}^{k - 1}} f \left( N \left( m + \left( x_{a_1} - x_{a_2}, x_{a_2} - x_{a_3}, \dots, x_{a_{k - 1}} - x_{a_k} \right) \right) \right) = \int_{\mathbb{R}^{k - 1}} f(x) \mathrm{d} x$$
This is equivalent to our definition, by first defining $g \left( y_1, \dots, y_{k - 1} \right) = f \left( y_1, y_1 + y_2, \dots, y_1 + \cdots + y_{k - 1} \right)$ which is also smooth and compactly supported, and in one direction approximating $g$ from above and below by a linear combination of indicators of boxes, and in the other direction approximating an indicator function of a box from above and below by a smooth test function.

Theorem \ref{real-result} generalizes in a straightforward manner in the following way:

\begin{thm}\label{k-level}

For all $c > 0$, there exists a well spaced sequence $x_n$ with additive energy $E_{N, \gamma} \gg N^3$ and metric Poissonian $k$-level correlation for all $k \geq 2$, where $\gamma = \frac{\log N \left( \log \log N \right)^{1 + c}}{N}$.

\end{thm}

This follows in the same manner from a generalization of Theorem \ref{disc_pois}:

\begin{thm}\label{disc-k}

Under the conditions of Theorem \ref{disc_pois}, almost surely the sequence $x_n + z_n$ has metric Poissonian $k$-level correlation for all $k \geq 2$.

\end{thm}

The proof is very similar to the proof of Theorem \ref{disc_pois}, and so we will just sketch what the differences are and how to overcome them.

The main reason for studying higher level correlations is that the $k$-level correlations for all $k \geq 2$ essentially determine the local statistics entirely. More precisely, let $x_{(1)}^{N} \leq x_{(2)}^{N} \leq \cdots \leq x_{(N)}^{N}$ denote the first $N$ elements of the sequence $x_n$, ordered by size, and define
$$h \left( x, \left( x_n \right)_{n \geq 1}, N \right) = \frac{1}{N} \# \left\{ 1 \leq n < N \ : \ N \left( x_{(n)}^{N} - x_{(n + 1)}^{N} \right) \leq x \right\}$$
The \emph{level spacing distribution} (also called the \emph{gap distribution}, or the \emph{nearest neighbour spacing distribution}) is the limiting distribution $P(s)$ (if it exists), which satisfies
$$\lim_{N \to \infty} h \left( x, \left( x_n \right)_{n \geq 1}, N \right) = \intop_{0}^{x} P(s) \mathrm{d} s$$
The level spacing distribution contains the majority of local information at scale $\frac{1}{N}$ of the sequence $x_n$, and so it is naturally an object of interest. $N$ random points on the unit circle approximate a Poisson point process, and so in the random case the level spacing distribution is asymptotically a Poisson distribution, that is $P(s) = e^{- s}$. Accordingly, a sequence is said to have \emph{Poissonian level spacing distribution} if its level spacing distribution is $P(s) = e^{- s}$. It is well known that Poissonian $k$-level correlation for all $k \geq 2$ implies Poissonian level spacing distribution. Thus, \ref{k-level} and \ref{disc-k} immediately imply

\begin{cor}\label{arithmetic_cor}

For all $c > 0$, there exists a well spaced sequence $x_n$ with additive energy $E_{N, \gamma} \gg N^3$ and metric Poissonian level spacing distribution for all $k \geq 2$, where $\gamma = \frac{\log N \left( \log \log N \right)^{1 + c}}{N}$.

\end{cor}

\begin{cor}\label{random_cor}

Under the conditions of Theorem \ref{disc_pois}, almost surely the sequence $x_n + z_n$ has metric Poissonian level spacing distribution for all $k \geq 2$.

\end{cor}

\textbf{Acknowledgements}: The authors thank Zeev Rudnick for helpful comments, in particular pointing out Corollary \ref{arithmetic_cor}. The second author received funding from the European Research Council (ERC) under the European Union’s Horizon 2020 research and innovation programme (Grant agreement No. 786758).

\section{Proof of Theorem \ref{disc_pois}}

Define the random variable
$$X_{s, N} = \frac{1}{N} \sum_{1 \leq n \neq m \leq N} \mathbbm{1}_{s, N} \left( x_n + z_n - x_m - z_m \right)$$
By definition, we have Poissonian pair correlation if and only if almost surely
$$\lim_{N \to \infty} X_{s, N} = 2 s$$
It is enough to show this for each $s$ individually, as if this is true for all $s \in \mathbb{Q}$ then it is true for all $s$, and a countable intersection of sets of full measure is a set of full measure. For simplicity, from here on we will denote $X_N = X_{s, N}$. We will prove this via a second moment method. First, we compute the expectation of $X_N$.

$$\mathbb{E} \left[ X_N \right] = \frac{1}{N} \sum_{1 \leq m \neq n \leq N} \mathbb{E} \left[ \mathbbm{1}_{s, N} \left( x_n + z_n - x_m - z_m \right) \right]$$
For $x \in \mathbb{R}$, denote
$$\rho_N (x) = \frac{1}{N} \sum_{n = 1}^{N} \frac{1}{2 g(n)} \mathbbm{1} \left( x \in \left[ x_n - g(n), x_n + g(n) \right] \right) = \frac{1}{N} \sum_{n = 1}^{N} \frac{1}{2 g(n)} \mathbbm{1} \left( x_n \in \left[ x - g(n), x + g(n) \right] \right)$$
One should think of $\rho_N (x)$ as the density function of the probability measure which corresponds to balls around $x_n$ of radius $g(n)$, each with a weight of $1/N$, that is for every function $f$ we have
$$\frac{1}{N} \sum_{n = 1}^{N} \mathbb{E} \left[ f \left( x_n + z_n \right) \right] = \intop_{0}^{1} f \left( x \right) \rho_N (x) \mathrm{d} x$$
Now, we define
$$h_{s, N} (x) = \sum_{n = 1}^{N} \mathbb{E} \left[ \mathbbm{1}_{s, N} \left( x_n + z_n - x \right) \right]$$
By using linearity of expectation,
$$\mathbb{E} \left[ X_N \right] - \frac{1}{N} \sum_{m = 1}^{N} \mathbb{E} \left[ h_{s, N} \left( x_m + z_m \right) \right] = \frac{1}{N} \sum_{m = 1}^{N} \mathbb{E} \left[ \mathbbm{1}_{s, N} \left( z_m - z_{m}' \right) \right]$$
where $z_{1}', z_{2}', \dots$ is an independent copy of $z_1, z_2, \dots$. Clearly,
$$\mathbb{E} \left[ \mathbbm{1}_{s, N} \left( z_m - z_{m}' \right) \right] \leq \frac{s}{N g(m)} \leq \frac{s}{N g(N)}$$
and so
$$\frac{1}{N} \sum_{m = 1}^{N} \mathbb{E} \left[ \mathbbm{1}_{s, N} \left( z_m - z_{m}' \right) \right] \leq \frac{s}{N g(N)} = o_N \left( 1 \right)$$
Therefore,
$$\mathbb{E} \left[ X_N \right] = \intop_{0}^{1} h_{s, N} (x) \rho_N (x) \mathrm{d} x + o_N \left( 1 \right)$$
What we need to show is that $\intop_{0}^{1} h_{s, N} (x) \rho_N (x) \mathrm{d} x$ converges to $2 s$ as $N \to \infty$. We will do this by showing that uniformly in $x$, we have $\rho_N (x) \to 1, \ h_{s, N} (x) \to 2 s$. Let us start with $\rho_N (x)$. This will follow from the fact that the discrepancy of $x_n$ is much smaller than $g(n)$: essentially, $\rho_N(x)$ counts the number of points $x_n$ in a certain interval which is of a length much larger than the discrepancy, and so we have the correct asymptotic for the number of $x_n$'s in that interval. There is a slight complication, due to the fact that the interval is not fixed. To remedy this, we will divide our sum into intervals of the form $L \leq n \leq L + k(L)$ for some $L, k(L)$, where we need to ensure simultaneously that $g$ is approximately the same at both endpoints of the range of summation (which allows us to treat the interval $\left[ x - g(n), x + g(n) \right]$ as fixed), and that the range of summation is long enough so that we can use our discrepancy bounds.

Specifically, we will take $k(L)$ small enough such that $\frac{g \left( L + k(L) \right)}{g(L)} \to_{L \to \infty} 1$. Then,
$$\sum_{L < n \leq L + k(L)} \frac{1}{2 g(n)} \mathbbm{1} \left( x_n \in \left[ x - g(n), x + g(n) \right] \right) \sim \frac{1}{2 g(L)} \sum_{L < n \leq L + k(L)} \mathbbm{1} \left( x_n \in \left[ x - g(n), x + g(n) \right] \right)$$
Because $g$ is monotonically decreasing,
$$\mathbbm{1} \left( x_n \in \left[ x - g \left( L + k(L) \right), x + g \left( L + k(L) \right) \right] \right) \leq \mathbbm{1} \left( x_n \in \left[ x - g(n), x + g(n) \right] \right) \leq \mathbbm{1} \left( x_n \in \left[ x - g(L), x + g(L) \right] \right)$$
Letting $c$ be either $g(L)$ or $g \left( L + k(L) \right)$, we will show that
$$\frac{1}{2 g(L)} \sum_{L < n \leq L + k(L)} \mathbbm{1} \left( x_n \in \left[ x - c, x + c \right] \right) \sim \frac{c k(L)}{g(L)} \sim k(L)$$
which by our argument in the previous paragraph ($c = g(L)$ functions as an upper bound, and $c = g \left( L + k(L) \right)$ functions as a lower bound), implies that
$$\sum_{L < n \leq L + k(L)} \frac{1}{2 g(n)} \mathbbm{1} \left( x_n \in \left[ x - g(n), x + g(n) \right] \right) \sim k(L)$$
which shows that $\rho_N (x) \to 1$ uniformly in $x$, as required. All we have left is to compute
$$\frac{1}{2 g(L)} \sum_{L < n \leq L + k(L)} \mathbbm{1} \left( x_n \in \left[ x - c, x + c \right] \right) = \frac{1}{2 g(L)} \left( \sum_{1 \leq n \leq L + k(L)} \mathbbm{1} \left( x_n \in \left[ x - c, x + c \right] \right) - \sum_{1 \leq n \leq L} \mathbbm{1} \left( x_n \in \left[ x - c, x + c \right] \right) \right)$$
By the definition of discrepancy,
$$\abs{\sum_{1 \leq n \leq L + k(L)} \mathbbm{1} \left( x_n \in \left[ x - c, x + c \right] \right) - 2 c \left( L + k(L) \right)} \leq \left( L + k(L) \right) D_{L + k(L)}$$
and
$$\abs{\sum_{1 \leq n \leq L} \mathbbm{1} \left( x_n \in \left[ x - c, x + c \right] \right) - 2 c L} \leq L D_{L}$$
and so
$$\abs{\frac{1}{2 g(L)} \sum_{L < n \leq L + k(L)} \mathbbm{1} \left( x_n \in \left[ x - c, x + c \right] \right) - 2 c k(L)} \leq \frac{L D_L + \left( L + k(L) \right) D_{L + k(L)}}{g(L)}$$
In order for our sum to be asymptotically equal to $2 c k(L)$, we need to have
$$\lim_{L \to \infty} \frac{L D_L}{g(L) k(L)} = 0$$
Choosing $k(L) = f(L) \frac{d_L}{g(L)}$ where $f(L)$ is some slowly growing function, by our regularity condition on $g$ we get the desired result.

The proof that $h_{s, N} (x) \to 2 s$ is very similar. Notice that
$$\frac{2 s}{N} \sum_{1 \leq n \leq N} \frac{1}{2 g(n)} \mathbbm{1} \left( x_n \in \left[ x - g(n) + \frac{s}{N}, x + g(n) - \frac{s}{N} \right] \right) \leq h_{s, N} (x) \leq$$
$$\leq \frac{2 s}{N} \sum_{1 \leq n \leq N} \frac{1}{2 g(n)} \mathbbm{1} \left( x_n \in \left[ x - g(n) - \frac{s}{N}, x + g(n) + \frac{s}{N} \right] \right)$$
where the lower bound follows from the fact that if $x_n \in \left[ x - g(n) + \frac{s}{N}, x + g(n) - \frac{s}{N} \right]$ then $\mathbb{E} \left[ \mathbbm{1}_{s, N} \left( x_n + z_n - x \right) \right] = \frac{s}{N g(n)}$, and the upper bound follows from the fact that if $\mathbbm{1}_{s, N} \left( x - x_n - z_n \right)$ is positive with a non zero probability, then $\abs{x_n - x} \leq g(n) + \frac{s}{N}$. Now, following exactly the same proof as before with $g(n) \pm \frac{s}{N}$ instead of $g(n)$ we get the desired result. Note that this is the point in the proof where we use the assumption $N g(N) \to \infty$: this is precisely saying that the $\pm \frac{s}{N}$ is negligible.

Now, we bound the variance:
$$\mathrm{Var} \left( X_N \right) = \mathbb{E} \left[ X_{N}^2 \right] - \mathbb{E} \left[ X_N \right]^2 =$$
$$= \frac{1}{N^2} \sum_{\substack{1 \leq m_1, m_2, n_1, n_2 \leq N \\ m_1 \neq n_1, \ m_2 \neq n_2}} \mathbb{E} \left[ \mathbbm{1}_{s, N} \left( x_{n_1} + z_{n_1} - x_{m_1} - z_{m_1} \right) \mathbbm{1}_{s, N} \left( x_{n_2} + z_{n_2} - x_{m_2} - z_{m_2} \right) \right] - \mathbb{E} \left[ X_N \right]^2$$
Notice that if $m_1, n_1, m_2, n_2$ are distinct, the random variables $z_{n_1} - z_{m_1}$ and $z_{n_2} - z_{m_2}$ are independent, and therefore
$$\frac{1}{N^2} \sum_{\substack{1 \leq m_1, m_2, n_1, n_2 \leq N \\ m_1, m_2, n_1, n_2 \ \mathrm{distinct}}} \mathbb{E} \left[ \mathbbm{1}_{s, N} \left( x_{n_1} + z_{n_1} - x_{m_1} - z_{m_1} \right) \mathbbm{1}_{s, N} \left( x_{n_2} + z_{n_2} - x_{m_2} - z_{m_2} \right) \right] =$$
$$= \frac{1}{N^2} \sum_{\substack{1 \leq m_1, m_2,
n_1, n_2 \leq N \\ m_1, m_2, n_1, n_2 \ \mathrm{distinct}}} \mathbb{E} \left[ \mathbbm{1}_{s, N} \left( x_{n_1} + z_{n_1} - x_{m_1} - z_{m_1} \right) \right] \mathbb{E} \left[ \mathbbm{1}_{s, N} \left( x_{n_2} + z_{n_2} - x_{m_2} - z_{m_2} \right) \right] \leq \mathbb{E} \left[ X_N \right]^2$$
which means that
$$\mathrm{Var} \left( X_N \right) \leq \frac{1}{N^2} \sum_{1 \leq m_1, n_1, m_2, n_2 \leq N}' \mathbb{E} \left[ \mathbbm{1}_{s, N} \left( x_{n_1} + z_{n_1} - x_{m_1} - z_{m_1} \right) \mathbbm{1}_{s, N} \left( x_{n_2} + z_{n_2} - x_{m_2} - z_{m_2} \right) \right]$$
where $\sum'$ means that the summation is taken over quadruples $m_1, n_1, m_2, n_2$ such that $m_1 \neq n_1, \ m_2 \neq n_2$ and the sets $\left\{ m_1, n_1 \right\}, \ \left\{ m_2, n_2 \right\}$ have a nonempty intersection. If the sets have $2$ elements in common, then $m_1 = m_2, \ n_1 = n_2$ or $m_1 = n_2, \ m_2 = n_1$, and so the sum over those quadruples is simply
$$\frac{2 \mathbb{E} \left[ X_N \right]}{N}$$
If the sets have an intersection of size $1$, then the sum is
$$\frac{4}{N^2} \sum_{\substack{1 \leq m, n_1, n_2 \leq N \\ m \neq n_1 \neq n_2 \neq m}} \mathbb{E} \left[ \mathbbm{1}_{s, N} \left( x_{n_1} + z_{n_1} - x_m - z_m \right) \mathbbm{1}_{s, N} \left( x_{n_2} + z_{n_2} - x_m - z_m \right) \right]$$
Using the identity
$$\frac{1}{N} \sum_{m = 1}^{N} \mathbb{E} \left[ f \left( x_m + z_m \right) \right] = \intop_{0}^{1} f(x) \rho_N (x) \mathrm{d} x$$
a few times, we get that the variance is bounded by
$$4 N \left( \intop_{0}^{1} \mathbbm{1}_{s, N} \left( x \right) \rho_N (x) \mathrm{d} x \right)^2 \left( \intop_{0}^{1} \rho_N (x) \mathrm{d} x \right) \sim \frac{16 s}{N}$$

By Chebyshev's inequality we have
$$\mathbb{P} \left( | X_N - \mathbb{E} \left[ X_N \right] | \geq \frac{1}{N^{1/4}} \right) \ll \frac{1}{N^{1/2}}$$
and so by Borel-Cantelli, almost surely we have pair correlation along the subsequence $X_{N^4}$. However, a standard argument shows that this is sufficient to get pair correlation, for example see \cite[Lemma 3.1]{pair-lem}. This concludes the proof of the theorem.

Let us note an interesting fact about the proof: when computing the expectation, we showed that $\rho_N (x)$ converges uniformly to $1$, whereas at first glance it seems that we needed a much weaker result. However, as we have seen above, we usually have $h_{s, N} (x) \sim 2 s \rho_N (x)$ and therefore
$$\intop_{0}^{1} h_{s, N} (x) \rho_N (x) \mathrm{d} x \sim 2 s \intop_{0}^{1} \rho_N (x)^2 \mathrm{d} x$$
This integral converges to $2 s$ if and only if
$$\intop_{0}^{1} \rho_N (x)^2 \mathrm{d} x \to 1$$
Recall that
$$\intop_{0}^{1} \rho_N (x) \mathrm{d} x = 1$$
and so we have the required convergence if and only if $\rho_N (x)$ converges to the constant $1$ in $L^2 \left[ 0, 1 \right]$. The variance bound is actually much more flexible than the expectation: all we need for the Chebyshev/Borel-Cantelli argument to work is
$$\mathrm{Var} \left( X_N \right) \ll \frac{1}{N^{\varepsilon}}$$
and in fact this can be further weakened to
$$\mathrm{Var} \left( X_N \right) \ll \frac{1}{\left( \log N \right)^{1 + \varepsilon}}$$
Therefore, if we can show that $\rho_N (x)$ does not converge to $1$ in $L^2$, which is reasonable to expect if $g(N)$ is smaller than the discrepancy, we can show that almost surely the pair correlation is too large to be Poissonian, and in particular we get that almost surely we do not have Poissonian pair correlation. In the next section, we shall partially follow this through for the sequence $x_n = \alpha n$ for almost all $\alpha$, but this is still far from a full converse. We think this is interesting, and a satisfactory converse theorem to Theorem \ref{disc_pois} would confirm this connection between discrepancy and pair correlation.

\section{Proof of Theorem \ref{converse}}

We set up the same machinery as in the previous section. We show that almost surely,
$$\limsup_{N \to \infty} X_{s, N} > 2 s$$
which will imply lack of Poissonian pair correlation. Again, we do this via a second moment method. As in the previous section, we have
$$\mathbb{E} \left[ X_N \right] = \intop_{0}^{1} h_{s, N} (x) \rho_N (x) \mathrm{d} x + o_N (1)$$
As we have sketched in the previous section, $h_{s, N} \sim 2 s \rho_N (x)$ and so if we show that $\rho_N$ does not converge to $1$ in $L^2$ we will have
$$\limsup_{N \to \infty} \mathbb{E} \left[ X_N \right] > 2 s$$
By Khintchine's theorem, for almost all $\alpha$, there exist infinitely many coprime $p, q$ such that
$$\abs{\alpha - \frac{p}{q}} \leq \frac{1}{q^2 \log q \log \log q}$$
Let $N = q \left( \log q \right)^{1/2} \left( \log \log q \right)^{1/3}$. From now on in this section, we will work only with $N$ along this specific subsequence. Then, for $1 \leq n \leq N$
$$\norm{n \alpha - n \frac{p}{q}} \leq \frac{N}{N^2 \left( \log \log q \right)^{1/3}} = o \left( \frac{1}{N} \right)$$
and therefore in studying $X_N$ we can replace the points $n \alpha$ with $n \frac{p}{q}$, as local statistics at scale $\frac{1}{N}$ are unchanged. Assume for contradiction that $\rho_N (x)$ converges to $1$ in $L^2$, and that $\rho_{N / 2}$ converges to $1$ in $L^2$ as well. Then, so does $2 \rho_N (x) - \rho_{N / 2} (x)$, which means that
$$\frac{2}{N} \sum_{n = N / 2}^{N} \frac{1}{2 g(n)} \mathbbm{1} \left( x \in \left[ n \frac{p}{q} - g(n), n \frac{p}{q} + g(n) \right] \right)$$
converges to $1$ in $L^2$. However, notice that the measure of the support of this function tends to $0$, because its support is a union of intervals centered around at most $q \leq \frac{N}{\left( \log N \right)^{1/2} \left( \log \log N \right)^{1/4}}$ distinct points, and each distinct point contributes to the set an interval of measure at most $\frac{\left( \log N \right)^{1/2}}{N}$. Therefore, $2 \rho_N - \rho_{N / 2}$ converges in measure to $0$, and in particular does not converge to $1$ in $L^2$.

Now we bound the variance. As in the previous section, it is bounded by
$$4 N \left( \intop_{0}^{1} \mathbbm{1}_{s, N} \left( x \right) \rho_N (x) \mathrm{d} x \right)^2$$
We will now give a slightly coarse pointwise bound on $\rho_N (x)$: we show that uniformly in $0 \leq x \leq 1$ we have
$$\rho_N (x) \ll \log N$$
We do this by splitting the sum into intervals of the form $L \leq n \leq 2 L$. We now need to bound
$$\frac{1}{g (2 L)} \sum_{n = 1}^{L} \mathbb{1} \left( x \in \left[ n \frac{p}{q} - g \left( L \right), n \frac{p}{q} + g \left( L \right) \right] \right)$$
Clearly, this is at most
$$\frac{g \left( L \right)}{g \left( 2 L \right)} \left( 1 + \frac{L}{q} \right) \ll 1 + \frac{L}{q}$$
Summing these $\log N$ intervals we get the desired result. This implies that we have
$$\mathrm{Var} \left( X_N \right) \ll \frac{\left( \log N \right)^2}{N}$$

Now applying Chebyshev's inequality and Borel-Cantelli along a subsequence shows that almost surely,
$$\limsup_{N \to \infty} X_{s, N} > 2 s$$
as required.
\section{Sketch of proof of Theorem \ref{disc-k}}

It is enough to show that for each $k$ individually, almost surely we have $k$-level Poissonian correlation, because a countable intersection of sets of full measure is a set of full measure. Define
$$X_{s, N} = \frac{1}{N} \sum_{a \in \mathcal{X}_k} \mathbbm{1}_{s, N} \left( x_{a_1} - x_{a_2}, x_{a_1} - x_{a_3}, \dots, x_{a_1} - x_{a_k} \right)$$
As before, we have Poissonian pair correlation if and only if
$$\lim_{N \to \infty} X_{s, N} = \left( s_2 - s_1 \right) \left( s_4 - s_3 \right) \cdots \left( s_{2 k - 2} - s_{2 k - 3} \right)$$
It is sufficient to prove this for each $s$ individually, and once again we use the second moment method. The expectation now ends up being asymptotically equivalent to
$$\left( s_2 - s_1 \right) \left( s_4 - s_3 \right) \cdots \left( s_{2 k - 2} - s_{2 k - 3} \right) \intop_{0}^{1} \rho_N (x)^k \mathrm{d} x$$
and and the asymptotic still follows from the fact that $\rho_N$ uniformly converges to $1$. As for the variance, the way to take care of it is to split $\mathbb{E} \left[ X_{s, N}^2 \right]$ into sums of the form
$$\sum_{\substack{a, b \in \mathcal{X}_k \\ \abs{a \cap b} = \ell}} \cdots$$
that is, the summation is over pairs of sequences $a, b$ with exactly $\ell$ common elements, exactly like what we did above. The contribution of the sum for $\ell = 0$ is at most $\mathbb{E} \left[ X_{s, N} \right]^2$ in the same manner as we did for pair correlation, and the contribution for each $\ell \geq 1$ is easily seen to be (by the same method) $\ll \frac{1}{N}$, and using Chebyshev's inequality and Borel-Cantelli proves the theorem.


\begin{thebibliography}{99}
\bibitem{2009.08184}
C. Aistleitner, D. El-Baz and M. Munsch, \textit{A pair correlation problem, and counting lattice points with the zeta function}, Geom. Funct. Anal. (2021), \url{https://doi.org/10.1007/s00039-021-00564-6}

\bibitem{pair->equi1}
C. Aistleitner, T. Lachmann and F. Pausinger, \textit{Pair correlations and equidistribution}, J. Number Theory 182, 206-220 (2018)

\bibitem{disc-ener}
C. Aistleitner, G. Larcher, \textit{Additive Energy and Irregularities of Distribution}, Uniform distribution theory, Volume 12, Issue 1, pp. 99-107 (2017).

\bibitem{bourgain}
C. Aistleitner, G. Larcher, and M. Lewko, \textit{Additive energy and the Hausdorff dimension of the exceptional set in metric pair correlation problems}, Isr. J. Math. 222, 463–485 (2017).

\bibitem{beck}
J. Beck, \textit{Probabilistic diophantine approximation, I. Kronecker sequences}, Ann. of Math. 140 (1994) 451–502.

\bibitem{khintchine}
A. Khintchine, \textit{Ein Satz \"{u}ber Kettenbr\"{u}che mit arithmetischen Anwendungen}, Math. Z. 18 (1923) 289–306.

\bibitem{book}
L. Kuipers, H. Niederrieter, \textit{Uniform distribution of sequences}, Pure and Applied Mathematics. Wiley-Interscience [John Wiley \& Sons], New York-London-Sydney, 1974.

\bibitem{lachmann}
T. Lachmann, N. Technau, \textit{On exceptional sets in the metric Poissonian pair correlations problem}. Monatsh Math 189, 137–156 (2019).

\bibitem{pair->equi2}
G. Larcher, S. Grepstad, \textit{On pair correlation and discrepancy}, Arch. Math., vol. 109 (2), 143–149, 2017.

\bibitem{bourgain+}
G. Larcher, W. Stockinger, \textit{Pair correlation of sequences $\left( a_n \alpha \right)_{n \in \mathbb{N}}$ with maximal additive energy}. Mathematical Proceedings of the Cambridge Philosophical Society, 168(2), 287-293, (2020).

\bibitem{neg}
G. Larcher, W. Stockinger, \textit{Some negative results related to Poissonian pair correlation problems}, Discrete Mathematics, Vol. 343, Nr. 2, pp. 18, (2020).

\bibitem{very small theta}
C. Lutsko, A. Sourmelidis, N. Technau, \textit{Pair Correlation of the Fractional Parts of $\alpha n^{\theta}$}, (2021), arXiv:2106.09800.

\bibitem{lowdisc}
W. Schmidt, \textit{Irregularities of distribution, VII}. Acta Arith. 21, 45–50, 1972

\bibitem{stein}
S. Steinerberger, \textit{Poissonian pair correlation and discrepancy}, Indagationes Mathematicae, Volume 29, Issue 5, 1167-1178, (2018).

\bibitem{pair-lem}
N. Technau, Z. Rudnick, \textit{The metric theory of the pair correlation function of real-valued lacunary sequences}. Illinois J. Math. Np. 64 (2020): 583-594.

\bibitem{small theta}
N. Technau, Z. Rudnick, \textit{The metric theory of the pair correlation function for small non-integer powers}. (2021), arXiv:2107.07092

\end{thebibliography}
\end{document}